\documentclass[12pt]{article}%
\usepackage{graphicx}
\usepackage[intlimits]{amsmath}
\usepackage{latexsym}
\usepackage{amsfonts}
\usepackage{amssymb}%
\makeatletter
\newfont{\footsc}{cmcsc10 at 8truept}
\newfont{\footbf}{cmbx10 at 8truept}
\newfont{\footrm}{cmr10 at 10truept}
\newtheorem{theorem}{Theorem}

\newtheorem{corollary}[theorem]{Corollary}

\newtheorem{lemma}[theorem]{Lemma}

\newtheorem{proposition}[theorem]{Proposition}

\newenvironment{proof}[1][Proof]{\noindent{\textbf {#1}  }}  {\hfill$\Box$\bigskip}

\begin{document}

\title{Eigenvalues and forbidden subgraphs I}
\author{Vladimir Nikiforov\\Department of Mathematical Sciences, University of Memphis, \\Memphis TN 38152, USA}
\maketitle

\begin{abstract}
Suppose a graph $G$ have $n$ vertices, $m$ edges, and $t$ triangles. Letting
$\lambda_{n}\left(  G\right)  $ be the largest eigenvalue of the Laplacian of
$G$ and $\mu_{n}\left(  G\right)  $ be the smallest eigenvalue of its
adjacency matrix, we prove that
\begin{align*}
\lambda_{n}\left(  G\right)   &  \geq\frac{2m^{2}-3nt}{m\left(  n^{2}%
-2m\right)  },\\
\mu_{n}\left(  G\right)   &  \leq\frac{3n^{3}t-4m^{3}}{nm\left(
n^{2}-2m\right)  },
\end{align*}
with equality if and only if $G$ is a regular complete multipartite graph.

Moreover, if $G$ is $K_{r+1}$-free, then
\[
\lambda_{n}\left(  G\right)  \geq\frac{2mn}{\left(  r-1\right)  \left(
n^{2}-2m\right)  }%
\]
with equality if and only if $G$ is a regular complete $r$-partite graph.


\end{abstract}

\section{Introduction}

Our notation is standard (e.g., see \cite{Bol98} and \cite{CDS80}); in
particular, $G\left(  n\right)  $ stands for a graph of order $n,$ and
$G\left(  n,m\right)  $ stands for a graph of order $n$ and size $m$. We write
$t\left(  G\right)  $ for the number of triangles of a graph $G,$ $A\left(
G\right)  $ for its adjacency matrix, and $D\left(  G\right)  $ for the
diagonal matrix of its degree sequence. The Laplacian of $G$ is defined as
$L\left(  G\right)  =D\left(  G\right)  -A\left(  G\right)  .$ Given a graph
$G=G\left(  n\right)  ,$ the eigenvalues of $A\left(  G\right)  $ are $\mu
_{1}\left(  G\right)  \geq...\geq\mu_{n}\left(  G\right)  $ and the
eigenvalues of $L\left(  G\right)  $ are $0=\lambda_{1}\left(  G\right)
\leq...\leq\lambda_{n}\left(  G\right)  .$

In this note we study how $\lambda_{n}\left(  G\right)  $ and $\mu_{n}\left(
G\right)  $ depend on the number of certain subgraphs of $G.$ In \cite{Nik06}
we showed that if $r\geq2$ and $G$ is a $K_{r+1}$-free graph with $n$ vertices
and $m$ edges, then
\begin{equation}
\mu_{n}\left(  G\right)  <-\frac{2}{r}\left(  \frac{2m}{n^{2}}\right)  ^{r}n.
\label{in1}%
\end{equation}

Here we prove a similar inequality for $\lambda_{n}\left(  G\right)  $.

\begin{theorem}
\label{th1} If $r\geq2$ and $G=G\left(  n,m\right)  $ is a $K_{r+1}$-free
graph, then
\begin{equation}
\lambda_{n}\left(  G\right)  \geq\frac{2mn}{\left(  r-1\right)  \left(
n^{2}-2m\right)  } \label{mbound}%
\end{equation}
with equality if and only if $G$ is a regular complete $r$-partite graph.
\end{theorem}

We deduce Theorem \ref{th1} from more general results.

\begin{theorem}
\label{th2}If $G=G\left(  n,m\right)  ,$ then%
\begin{equation}
6nt\left(  G\right)  \geq\left(  n+\lambda_{n}\left(  G\right)  \right)
\sum_{u\in V\left(  G\right)  }d^{2}\left(  u\right)  -2nm\lambda_{n}\left(
G\right)  \label{l1}%
\end{equation}
with equality if and only if $G$ is a complete multipartite graph, and
\begin{equation}
\lambda_{n}\left(  G\right)  \geq\frac{2m^{2}-3nt\left(  G\right)  }{m\left(
n^{2}-2m\right)  }n \label{l2}%
\end{equation}
with equality if and only if $G$ is a regular complete multipartite graph.
\end{theorem}

Inequality (\ref{l2}) suggests a similar inequality for $\mu_{n}\left(
G\right)  .$

\begin{theorem}
\label{th3}If $G=G\left(  n,m\right)  ,$ then%
\begin{equation}
\mu_{n}\left(  G\right)  \leq\frac{3n^{3}t\left(  G\right)  -4m^{3}}{nm\left(
n^{2}-2m\right)  } \label{in2}%
\end{equation}
with equality if and only if $G$ is a regular complete multipartite graph.
\end{theorem}

From Theorem \ref{th3} we effortlessly deduce results complementary to results
of Serre, Li, and Cioab\u{a} (e.g., see \cite{Cio06} and its references). Note
that these authors study regular graphs of fixed degree and large order; in
contrast, our results are meaningful for any graph $G=G\left(  n\right)  $
with average degree $\gg n^{1/2}.$ Here is an immediate consequence of Theorem
\ref{th3}.

\begin{corollary}
\label{th4} If $0\leq\varepsilon\leq1$ and $G=G\left(  n,m\right)  $ is a
graph with $t\left(  G\right)  \leq\varepsilon\left(  m/n\right)  ^{3}$, then%
\[
\mu_{n}\left(  G\right)  \leq-\left(  1-\varepsilon\right)  \frac{4m^{2}%
}{n^{3}}.
\]

\end{corollary}

In other words, graphs with small $\left\vert \mu_{n}\left(  G\right)
\right\vert $ abound in triangles. Likewise, graphs with small $\left\vert
\mu_{n}\left(  G\right)  \right\vert $ have cycles of all lengths up to
$O\left(  m^{2}/n^{3}\right)  $.

\begin{corollary}
\label{th5} If $3\leq r\leq n/2$ and a graph $G=G\left(  n,m\right)  $
contains no cycle of length $r$, then
\[
\mu_{n}\left(  G\right)  \leq-\frac{4m^{2}}{n^{3}}+2\left(  r-3\right)  .
\]

\end{corollary}

\section{Proofs}

For any vertex $u,$ write $\Gamma\left(  u\right)  $ for the set of its
neighbors, $t\left(  u\right)  $ for the number of triangles containing it,
and $t^{\prime}\left(  u\right)  $ for $e\left(  V\left(  G\right)
\backslash\Gamma\left(  u\right)  \right)  .$\bigskip

\begin{proof}
[Proof of Theorem \ref{th2}]It is known (e.g., see \cite{BoNi04}) that for any
partition $V\left(  G\right)  =V_{1}\cup V_{2},$
\begin{equation}
\lambda_{n}\left(  G\right)  \geq\frac{e\left(  V_{1},V_{2}\right)
n}{\left\vert V_{1}\right\vert \left\vert V_{2}\right\vert }. \label{mainin}%
\end{equation}
Therefore, for every $u\in V\left(  G\right)  $ and partition $V_{1}%
=\Gamma\left(  u\right)  $, $V_{2}=V\left(  G\right)  \backslash\Gamma\left(
u\right)  ,$%
\begin{equation}
\lambda_{n}\left(  G\right)  d\left(  u\right)  \left(  n-d\left(  u\right)
\right)  \geq ne\left(  \Gamma\left(  u\right)  ,V\left(  G\right)
\backslash\Gamma\left(  u\right)  \right)  . \label{in4}%
\end{equation}
In view of%
\[
\sum_{v\in\Gamma\left(  u\right)  }d\left(  v\right)  =e\left(  \Gamma\left(
u\right)  ,V\left(  G\right)  \backslash\Gamma\left(  u\right)  \right)
+2t\left(  u\right)  ,
\]
summing (\ref{in4}) for all $u\in V\left(  G\right)  ,$ we find that
\begin{align*}
\lambda_{n}\left(  G\right)  \sum_{u\in V\left(  G\right)  }d\left(  u\right)
\left(  n-d\left(  u\right)  \right)   &  \geq n\sum_{u\in V\left(  G\right)
}\sum_{v\in\Gamma\left(  u\right)  }d\left(  v\right)  -2n\sum_{u\in V\left(
G\right)  }2t\left(  u\right) \\
&  =n\sum_{u\in V\left(  G\right)  }d^{2}\left(  u\right)  -6nt\left(
G\right)  ,
\end{align*}
and (\ref{l1}) follows. Now (\ref{l2}) follows from (\ref{l1}), in view of
$\sum_{u\in V\left(  G\right)  }d^{2}\left(  u\right)  \geq4m^{2}/n.$

If equality holds in (\ref{l1}), then equality holds in (\ref{in4}) for every
$u\in V\left(  G\right)  .$ Note that if equality holds in (\ref{mainin}),
then all vertices from $V_{2}$ are connected to the same number of vertices
from $V_{1}$ (for a detailed proof of this result see, e.g., \cite{BoNi06}.)
In our selection of $V_{1}$ and $V_{2}$ the vertex $u$ is joined to all
vertices from $\Gamma\left(  u\right)  ,$ hence all vertices from $V\left(
G\right)  \backslash\Gamma\left(  u\right)  $ are joined to all vertices from
$\Gamma\left(  u\right)  .$ Consequently, $G$ contains no induced subgraph of
order $3$ with exactly one edge; hence, $G$ is complete multipartite.

If equality holds in (\ref{l2}), then $G$ is a complete multipartite graph; as
$\sum_{u\in V\left(  G\right)  }d^{2}\left(  u\right)  =4m^{2}/n,$ $G$ is regular.
\end{proof}

\begin{proof}
[\textbf{Proof of Theorem \ref{th1}}]Since $G$ is $K_{r+1}$-free,
$\Gamma\left(  u\right)  $ induces a $K_{r}$-free graph for every $u\in
V\left(  G\right)  .$ Thus, Tur\'{a}n's theorem implies that
\[
t\left(  u\right)  \leq\frac{r-2}{2\left(  r-1\right)  }d^{2}\left(  u\right)
.
\]
Summing this inequality for all $u\in V\left(  G\right)  ,$ we obtain%
\[
6t\left(  G\right)  \leq\frac{r-2}{r-1}\sum_{u\in V\left(  G\right)  }%
d^{2}\left(  u\right)  .
\]
This, in view of (\ref{l1}), implies that
\[
n\frac{r-2}{r-1}\sum_{u\in V\left(  G\right)  }d^{2}\left(  u\right)
\geq\left(  n+\lambda_{n}\left(  G\right)  \right)  \sum_{u\in V\left(
G\right)  }d^{2}\left(  u\right)  -2nm\lambda_{n}\left(  G\right)  .
\]
Using $\sum_{u\in V\left(  G\right)  }d^{2}\left(  u\right)  \geq4m^{2}/n,$
the result follows after simple algebra.

If equality holds in (\ref{mbound}), then equality holds in (\ref{l1}),
implying that $G$ is a complete multipartite graph. The condition for equality
in Tur\'{a}n's theorem implies that the neighborhood of every vertex is a
complete $\left(  r-1\right)  $-partite graph, thus, $G$ is $r$-partite.
Finally, we have $\sum_{u\in V\left(  G\right)  }d^{2}\left(  u\right)
=4m^{2}/n,$ so $G$ is regular, completing the proof.
\end{proof}

\subsection{Proof of Theorem \ref{th3}}

To prove Theorem \ref{th3}, we need two propositions and a lemma.

\begin{proposition}
\label{pro1} For every graph $G=G\left(  n,m\right)  ,$%
\[
2\sum_{u\in V\left(  G\right)  }d\left(  u\right)  \left(  t^{\prime}\left(
u\right)  -t\left(  u\right)  \right)  =4m^{2}-4\sum_{uv\in E\left(  G\right)
}d\left(  u\right)  d\left(  v\right)
\]

\end{proposition}

\begin{proof}
For every $u\in V\left(  G\right)  ,$ we have%
\begin{align*}
2t\left(  u\right)   &  =\sum_{v\in\Gamma\left(  u\right)  }d\left(  v\right)
-e\left(  V_{1},V_{2}\right) \\
2t^{\prime}\left(  u\right)   &  =\sum_{v\in V\left(  G\right)  \backslash
\Gamma\left(  u\right)  }d\left(  v\right)  -e\left(  V_{1},V_{2}\right)
=2m-e\left(  V_{1},V_{2}\right)  -\sum_{v\in\Gamma\left(  u\right)  }d\left(
v\right)  .
\end{align*}
Hence,%
\[
2\left(  t^{\prime}\left(  u\right)  -t\left(  u\right)  \right)  d\left(
u\right)  =2md\left(  u\right)  -2d\left(  u\right)  \sum_{v\in\Gamma\left(
u\right)  }d\left(  v\right)  ;
\]
summing this equality for all $u\in V\left(  G\right)  ,$ we obtain the
required equality$.$
\end{proof}

\bigskip

\begin{proposition}
\label{pro2} For every graph $G=G\left(  n,m\right)  $%
\[
2\sum_{uv\in E\left(  G\right)  }d\left(  u\right)  d\left(  v\right)
\geq4m^{2}+\sum_{u\in V\left(  G\right)  }d^{3}\left(  u\right)
-nd^{2}\left(  u\right)
\]

\end{proposition}

\begin{proof}
We have
\begin{align*}
2\sum_{uv\in E\left(  G\right)  }d\left(  u\right)  d\left(  v\right)   &
=\sum_{uv\in E\left(  G\right)  }d^{2}\left(  u\right)  +d^{2}\left(
v\right)  -\left(  d\left(  u\right)  -d\left(  v\right)  \right)  ^{2}\\
&  \geq\sum_{u\in V\left(  G\right)  }d^{3}\left(  u\right)  -\frac{1}{2}%
\sum_{u\in V\left(  G\right)  }\sum_{v\in V\left(  G\right)  }\left(  d\left(
u\right)  -d\left(  v\right)  \right)  ^{2}\\
&  =\sum_{u\in V\left(  G\right)  }d^{3}\left(  u\right)  -\frac{1}{2}%
\sum_{u\in V\left(  G\right)  }\sum_{v\in V\left(  G\right)  }d^{2}\left(
u\right)  +d^{2}\left(  v\right)  -2d\left(  u\right)  d\left(  v\right) \\
&  =\sum_{u\in V\left(  G\right)  }d^{3}\left(  u\right)  -n\sum_{u\in
V\left(  G\right)  }d^{2}\left(  u\right)  +4m^{2},
\end{align*}
completing the proof.
\end{proof}

\bigskip

\begin{lemma}
\label{le1}Let $0\leq x_{1}\leq\ldots\leq x_{n}\leq1$ be reals with
$x_{1}+\ldots+x_{n}=ns.$ Then,
\begin{equation}
\sum_{i=1}^{n}2x_{i}^{3}-\left(  2+s\right)  x_{i}^{2}\geq ns^{3}-2ns^{2}.
\label{in3}%
\end{equation}

\end{lemma}

\begin{proof}
Setting $\varphi\left(  x\right)  =2x^{3}-\left(  2+s\right)  x^{2},$ we
routinely find that:

\emph{(i)} $\varphi\left(  x\right)  $ decreases for $0\leq x\leq\left(
s+2\right)  /3$ and increases for $\left(  s+2\right)  /3\leq x\leq1;$

\emph{(ii)} $\varphi\left(  x\right)  $ is concave for $0\leq x\leq\left(
s+2\right)  /6;$

\emph{(iii)} $\varphi\left(  x\right)  $ is convex for $\left(  s+2\right)
/6\leq x\leq1.$

Let $F\left(  x_{1},\ldots,x_{n}\right)  =\sum_{i=1}^{n}\varphi\left(
x_{i}\right)  $ and suppose $x_{1}\leq\ldots\leq x_{n}$ are such that
$F\left(  x_{1},\ldots,x_{n}\right)  $ is minimal, subject to the conditions
of the lemma; clearly, we may assume that $x_{1}>0.$

If $x_{1}\geq\left(  s+2\right)  /6$, \emph{(iii)} implies that $x_{1}%
=\ldots=x_{n}$ and the proof is completed. Assume $x_{1}<\left(  s+2\right)
/6;$ we shall show that this assumption leads to a contradiction. Note first
that $x_{2}\geq\left(  s+2\right)  /6$; otherwise, for sufficiently small
$\varepsilon>0,$ \emph{(ii)} implies that%
\[
F\left(  x_{1}-\varepsilon,x_{2}+\varepsilon,x_{3,}\ldots,x_{n}\right)
<F\left(  x_{1},x_{2},x_{3,}\ldots,x_{n}\right)  ,
\]
contradicting the choice of $x_{1},\ldots,x_{n}.$ Using \emph{(iii)} again, we
find that $x_{2}=\ldots=x_{n}.$ Now, setting $z=x_{2},$ we see that
$0<x_{1}<s$ and $s<z<ns/\left(  n-1\right)  ,$ and that the function%
\begin{align*}
f\left(  z\right)   &  =F\left(  x_{1},z,\ldots,z\right)  \\
&  =2\left(  ns-\left(  n-1\right)  z\right)  ^{3}+2\left(  n-1\right)
z^{3}-\left(  2+s\right)  \left(  \left(  n-1\right)  z-ns\right)
^{2}-\left(  2+s\right)  \left(  n-1\right)  z^{2}%
\end{align*}
has a local minimum in the interval
\begin{equation}
s<z<ns/\left(  n-1\right)  .\label{intz}%
\end{equation}
We have%
\begin{align*}
f^{\prime}\left(  z\right)   &  =-6\left(  n-1\right)  \left(  \left(
n-1\right)  z-ns\right)  ^{2}+6\left(  n-1\right)  z^{2}\\
&  -2\left(  n-1\right)  \left(  2+s\right)  \left(  \left(  n-1\right)
z-ns\right)  -2\left(  2+s\right)  \left(  n-1\right)  z\\
&  =-6\left(  n-1\right)  \left(  n\left(  \left(  n-2\right)  z-ns\right)
\left(  z-s\right)  +\frac{\left(  2+s\right)  }{3}n\left(  z-s\right)
\right)  \\
&  =-6n\left(  n-1\right)  \left(  z-s\right)  \left(  \left(  n-2\right)
z-\left(  n-2\right)  s-2s+\frac{\left(  2+s\right)  }{3}\right)  \\
&  =-6n\left(  n-1\right)  \left(  n-2\right)  \left(  z-s\right)  \left(
z-\left(  s+\frac{5s-2}{3\left(  n-2\right)  }\right)  \right)  .
\end{align*}
In view of (\ref{intz}), the local minimum of $f\left(  z\right)  $ must be
attained at
\[
z_{0}=s+\frac{5s-2}{3\left(  n-2\right)  },
\]
implying, in particular, that $z_{0}>s$. But since $f^{\prime}\left(
z\right)  >0$ for $s<z<z_{0}$ and $f^{\prime}\left(  z\right)  <0$ for
$z>z_{0},$ we see that $f\left(  z\right)  $ has a local maximum at $z_{0}.$
This contradiction completes the proof.
\end{proof}

\bigskip

\begin{proof}
[\textbf{Proof of Theorem \ref{th3}}]In \cite{BoNi04} it is proved that for
any partition $V\left(  G\right)  =V_{1}\cup V_{2},$%
\begin{equation}
\mu_{n}\left(  G\right)  \leq\frac{2e\left(  V_{1}\right)  }{\left\vert
V_{1}\right\vert }+\frac{2e\left(  V_{2}\right)  }{\left\vert V_{2}\right\vert
}-\frac{2m}{n}. \label{mainin1}%
\end{equation}

Hence, for every $u\in V\left(  G\right)  $ and partition $V_{1}=\Gamma\left(
u\right)  $, $V_{2}=V\left(  G\right)  \backslash\Gamma\left(  u\right)  ,$
\begin{equation}
\mu_{n}\left(  G\right)  \leq\frac{2e\left(  V_{1}\right)  }{d\left(
u\right)  }+\frac{2e\left(  V_{2}\right)  }{n-d\left(  u\right)  }-\frac
{2m}{n}=\frac{2t\left(  u\right)  }{d\left(  u\right)  }+\frac{2t^{\prime
}\left(  u\right)  }{n-d\left(  u\right)  }-\frac{2m}{n}, \label{in6}%
\end{equation}
and therefore,
\begin{align*}
\mu_{n}\left(  G\right)  \left(  n-d\left(  u\right)  \right)  d\left(
u\right)   &  \leq2t\left(  u\right)  \left(  n-d\left(  u\right)  \right)
+2t^{\prime}\left(  u\right)  d\left(  u\right)  -\frac{2m}{n}d\left(
u\right)  \left(  n-d\left(  u\right)  \right) \\
&  =2nt\left(  u\right)  +d\left(  u\right)  \left(  2t^{\prime}\left(
u\right)  -2t\left(  u\right)  \right)  -2md\left(  u\right)  +\frac{2m}%
{n}d^{2}\left(  u\right)  .
\end{align*}
Summing this inequality for all $u\in V\left(  G\right)  $, in view of
$\mu_{n}\left(  G\right)  \leq0$ and
\begin{equation}
\sum_{u\in V\left(  G\right)  }\left(  n-d\left(  u\right)  \right)  d\left(
u\right)  \leq\frac{2m}{n}\left(  n^{2}-2m\right)  , \label{in7}%
\end{equation}
we obtain%
\[
\mu_{n}\left(  G\right)  \frac{2m}{n}\left(  n^{2}-2m\right)  =6nt\left(
G\right)  +\sum_{u\in V\left(  G\right)  }d\left(  u\right)  \left(
2t^{\prime}\left(  u\right)  -2t\left(  u\right)  \right)  -4m^{2}+\frac
{2m}{n}\sum_{u\in V\left(  G\right)  }d^{2}\left(  u\right)  .
\]
Propositions \ref{pro1} and \ref{pro2} imply that
\begin{align*}
\mu_{n}\left(  G\right)  \frac{2m}{n}\left(  n^{2}-2m\right)   &
\leq6nt\left(  G\right)  -4\sum_{uv\in E\left(  G\right)  }d\left(  u\right)
d\left(  v\right)  +\frac{2m}{n}\sum_{u\in V\left(  G\right)  }d^{2}\left(
u\right) \\
&  \leq6nt\left(  G\right)  -8m^{2}-2\sum_{u\in V\left(  G\right)  }%
d^{3}\left(  u\right)  +\left(  \frac{2m}{n}+2n\right)  \sum_{u\in V\left(
G\right)  }d^{2}\left(  u\right)  .
\end{align*}
Assume for convenience that $V\left(  G\right)  =\left\{  1,\ldots,n\right\}
$ and $d\left(  1\right)  \leq\ldots\leq d\left(  n\right)  .$ Setting
$x_{i}=d\left(  i\right)  /n,$ $1\leq i\leq n,$ and $s=2m/n^{2},$ Lemma
\ref{le1} implies that
\[
-2\sum_{u\in V\left(  G\right)  }d^{3}\left(  u\right)  +\left(  \frac{2m}%
{n}+2n\right)  \sum_{u\in V\left(  G\right)  }d^{2}\left(  u\right)
\leq8m^{2}-\frac{8m^{3}}{n^{2}};
\]
therefore,%
\[
\mu_{n}\left(  G\right)  \leq\frac{6nt\left(  G\right)  -8m^{3}/n^{2}}{\left(
2m/n\right)  \left(  n^{2}-2m\right)  },
\]
and (\ref{in2}) follows.

If equality holds\ in (\ref{in2}), then equality holds in (\ref{in7}); thus,
$G$ is regular. Also, equality holds in (\ref{in6}) for every $u\in V\left(
G\right)  .$ Some algebra shows that for regular graphs inequality
(\ref{mainin1}) is equivalent to (\ref{mainin}); hence, as in the proof of
Theorem \ref{th2}, $G$ is a complete multipartite graph. The proof is completed.
\end{proof}

\bigskip

\begin{proof}
[\textbf{Proof of Corollary \ref{th5}}]Let $G$ have $m$ edges and $t$
triangles. If $G$ has no $C_{r}$ for some $r\leq n/2,$ then $e\left(
G\right)  \leq n^{2}/4$ (\cite{Bol78}, p. 150). Since the neighborhood of any
vertex $u$ has no path of order $r-1$, by a theorem of Erd\H{o}s and Gallai
\cite{ErGa59}, the neigborhood of $u$ induces at most $\left(  r-3\right)
d\left(  u\right)  /2$ edges, i.e., $2t\left(  u\right)  \leq\left(
r-3\right)  d\left(  u\right)  $. Summing over all vertices, we see that
$3t\leq\left(  r-3\right)  m.$ Hence, Theorem \ref{th3} implies that
\[
\mu_{n}\left(  G\right)  \leq\frac{3n^{3}t-4m^{3}}{nm\left(  n^{2}-2m\right)
}\leq\frac{n^{3}\left(  r-3\right)  m-4m^{3}}{nm\left(  n^{2}-2m\right)  }%
\leq-\frac{4m^{2}}{n^{3}}+2\left(  r-3\right)  ,
\]
completing the proof.\bigskip
\end{proof}

\textbf{Acknowledgment }Part of this research was completed while the author
was visiting the Institute for Mathematical Sciences, National University of
Singapore in 2006. The author is also indebted to B\'{e}la Bollob\'{a}s for
his kind support.


\begin{thebibliography}{9}                                                                                                %


\bibitem {Bol78}B. Bollob\'{a}s, \emph{Extremal Graph Theory}, Academic Press, 1978.

\bibitem {Bol98}B. Bollob\'{a}s, \emph{Modern Graph Theory,} Graduate Texts in
Mathematics, \textbf{184,} Springer-Verlag, New York (1998), xiv+394 pp.

\bibitem {BoNi04}B. Bollob\'{a}s, V. Nikiforov, Graphs and Hermitian matrices:
eigenvalue interlacing, \emph{Discrete Math.} \textbf{289 }(2004), 119-127.

\bibitem {BoNi06}B. Bollob\'{a}s, V. Nikiforov, Graphs and Hermitian matrices:
exact interlacing, sumbitted.

\bibitem {CDS80}D. Cvetkovi\'{c}, M. Doob, H. Sachs, \emph{Spectra of Graphs,}
VEB Deutscher Verlag der Wissenschaften, Berlin, 1980, 368 pp.

\bibitem {Cio06}S. M. Cioab\u{a}, On the extreme eigenvalues of regular
graphs. \emph{J. Combin. Theory Ser. B} 96 (2006), 367--373.

\bibitem {ErGa59}P. Erd\H{o}s, T. Gallai On maximal paths and circuits of
graphs. \emph{Acta Math. Acad. Sci. Hungar} \textbf{10} (1959) 337--356.

\bibitem {Nik06}V. Nikiforov, The smallest eigenvalue of $K_{r}$-free graphs,
\emph{Discrete Math.} 306 (2006), 612-616.
\end{thebibliography}
\end{document}